\documentclass[reqno]{amsart}
\usepackage{amsmath,amsfonts,amsthm,amsopn,color,amssymb,enumitem}
\usepackage{palatino}
\usepackage{graphicx}
\usepackage[english]{babel}
\usepackage{caption}
\usepackage{subcaption}
\usepackage[colorlinks=true]{hyperref}
\hypersetup{urlcolor=blue, citecolor=red, linkcolor=blue}
\usepackage[utf8]{inputenc}

\usepackage{cleveref}
\crefname{equation}{}{}
\crefname{figure}{}{}

\usepackage{esint}
\usepackage[title]{appendix}

\usepackage[backend=bibtex,style=numeric,sorting=nyt,firstinits=true,isbn=false,url=false,doi=false]{biblatex}

\newtheorem{theorem}{Theorem}[section]
\newtheorem{lemma}[theorem]{Lemma}
\newtheorem{corollary}[theorem]{Corollary}

\newtheorem{proposition}[theorem]{Proposition}

\theoremstyle{definition}

\theoremstyle{remark}
\newtheorem{remark}[theorem]{Remark}

\def\beq{\begin{equation}}
\def\eeq{\end{equation}}
\def\ba{\begin{array}}
\def\ea{\end{array}}

\def\R{\mathbb R}

\newcommand{\rmnote}[1]{}

\numberwithin{equation}{section}

\newenvironment{key words}{\textbf{Keywords}\mbox{  }}{ }

\begin{document}



\title[Modica estimates and curvature results for overdetermined $p$-Laplace problems]{Modica type estimates and curvature results for overdetermined $p$-Laplace problems}

\author{Yuanyuan Lian}
\address{(Y.~Lian)
Departamento de An\'alisis matem\'atico, Universidad de Granada,
Campus Fuentenueva, 18071 Granada, Spain} \email{lianyuanyuan.hthk@gmail.com; yuanyuanlian@correo.ugr.es}

\author{Jing Wu}
\address{J. W.: Departamento de Matem\'aticas,  Universidad Aut\'onoma de Madrid, Ciudad Universitaria de Cantoblanco, 28049 Madrid, Spain}
\email{jing.wu@uam.es; jingwulx@126.com}

\date{\today}
\thanks{2020 {\it Mathematics Subject Classification.}
35N25, 35B50, 35J92}

\maketitle
\begin{abstract}
In this paper we prove Modica type estimates for the following overdetermined $p$-Laplace problem
\begin{equation*}
  \begin{cases}
  \mathrm{div} \left(|\nabla u|^{p-2}\nabla u\right)+f(u) =0& \mbox{in $\Omega$, }\\
  u>0 &\mbox{in $\Omega$, }\\
  u=0 &\mbox{on $\partial\Omega$, }\\
  \partial_{\nu} u=-\kappa &\mbox{on $\partial\Omega$, }
  \end{cases}
\end{equation*}
where $1<p<+\infty$, $f\in C^1(\mathbb{R})$, $\Omega \subset \mathbb{R}^n$ ($n\geq 2$) is a $C^1$ domain (bounded or unbounded), $\nu$ is the exterior unit normal of $\partial \Omega$ and $\kappa\geq 0$ is a constant. Based on Modica type estimates, we obtain rigidity results for bounded solutions. In particular, we prove that if there exists a nonpositive primitive $F$ of $f$ satisfying $F(0)\geq -(p-1)\kappa^p / p$ (for $p>2$ we also assume that if $F(u_0)=0$, $F(u)=O(|u-u_0|^p)$ as $u\rightarrow u_0$), then either the mean curvature of $\partial \Omega$ is strictly negative or $\Omega$ is a half-space.
\end{abstract}

\section{Introduction and statement of the results}
\label{Section 1}

We consider a $C^1$ domain (bounded or unbounded) $\Omega \subset\mathbb{R}^{n}$ ($n\geq 2$) and a solution $u$ of
\begin{equation}\label{e.P}
  \begin{cases}
  \mathrm{div} \left(|\nabla u|^{p-2}\nabla u\right)+f(u) =0& \mbox{in $\Omega$, }\\
  u>0 &\mbox{in $\Omega$, }\\
  u=0 &\mbox{on $\partial\Omega$, }\\
  \partial_{\nu} u=-\kappa &\mbox{on $\partial\Omega$, }
  \end{cases}
\end{equation}
where $1<p<+\infty$, $f\in C^1(\mathbb{R})$, $\nu$ denotes the exterior unit normal of $\partial \Omega$ and $\kappa\geq 0$ is a constant. The coexistence of Dirichlet and Neumann boundary conditions renders this problem overdetermined. Such problems have long attracted interest due to their profound implications in connecting analytic properties of partial differential equations to the geometry of underlying domains. A landmark result is Serrin's celebrated work \cite{MR333220}, which established that for $p=2$ in problem \eqref{e.P}, the existence of a solution in a bounded domain implies that the bounded domain $\Omega$ is a ball and the solution $u$ is radially symmetric. Weinberger \cite{MR333221} gave an alternative proof via integral identities and the maximum principle for an auxiliary function with $f\equiv1$. Serrin used the moving plane method which is more flexible for a general equation and Weinberger's method relaxes boundary regularity requirements and imposes weaker assumptions on the solution's regularity, see \cite{MR1200301}. These rigidity results find extensive applications in physics \cite{Sirakov2002}.

\medskip

It is natural to ask whether Serrin's result can be established for more general quasilinear equations. For equations with $p$-Laplacian, Garofalo and Lewis \cite{MR980297} extended the analogous rigidity result for problem \eqref{e.P} with $f=1$ by employing Weinberger's approach. Colesanti \cite{MR1381289} considered problem \cref{e.P} with $f(u,|\nabla u|)$ in a bounded $C^2$ domain and by the using moving plane method, he obtained if $u$ has only one critical point in $\Omega$ then $\Omega$ is a ball and $u$ is radially symmetric. Because of the degeneracy of $p$-Laplacian at the critical points of solutions, the difficulty of using the moving plane method is that the comparison principles are not available in the same way as the Laplacian, see \cite{MR1648566}. For $1<p<2$, if $f$ is locally Lipschitz, Damascelli and Pacella \cite{MR1776351} obtained the rigidity by the moving plane method. For $1<p<+\infty$, under assumptions of the convexity of $\Omega$ and a certain control of $f$, Brock and Henrot \cite{MR1947461} furnished a proof via Steiner symmetrization. Recently, the rigidity for the full range $1<p<+\infty$ was obtained by Cao, Wei and Zhan \cite{cao202} in general bounded domains.

In addition, since Reilly \cite{MR645791} found the relation between the Laplacian and the geometrical concept of mean curvature, the rigidity result of overdetermined problems can be obtained by the Soap Bubble Theorem. Via this idea, Colasuonno and Ferrari \cite{MR4043775} got the symmetric result of solutions of problem \cref{e.P} with $f=1$ and $1<p<2$. Many other refinements and generalizations about the rigidity results with $p$-Laplacian have been proposed (see \cite{MR2366129, MR3127383, MR3264508, MR4294646, MR4715314}).

\medskip

However, all the above results are in bounded domains. For unbounded domains, Berestycki, Caffarelli and Nirenberg \cite{MR1470317} considered the overdetermined problem with the Laplacian operator (i.e., $p=2$ in \eqref{e.P}):
\begin{equation}\label{e.Laplacian}
  \begin{cases}
  \Delta u+f(u) =0& \mbox{in $\Omega$, }\\
  u>0 &\mbox{in $\Omega$, }\\
  u=0 &\mbox{on $\partial\Omega$, }\\
  \partial_{\nu} u=-\kappa&\mbox{on $\partial\Omega$, }
  \end{cases}
\end{equation}
where $f$ is a Lipschitz function, $\kappa\geq 0$ is some constant and $\Omega\subset \mathbb{R}^n$ is unbounded with a Lipschitz graph. They proposed the following:

\medskip

\textbf{BCN Conjecture}: If $\Omega$ is a smooth domain with its complement $\Omega^c$ connected, then the existence of a bounded positive solution to problem \cref{e.Laplacian} implies that $\Omega$ is either a ball, a half-space, a generalized cylinder $B^k \times \mathbb{R}^{n-k}$ ($B^k$ is a ball in $\mathbb{R}^k$), or the complement of one of them.

\medskip

Research shows that the BCN Conjecture is false and the first counterexample is given by Sicbaldi \cite{MR2592974} with the topology of a cylinder. Since then, many counterexamples have been built by bifurcation arguments, see \cite{MR4484836, MR2854185,MR4649188} (topology of a cylinder),  \cite{MR3417183} (topology of a half-space) and \cite{MR4046014,dai2024} (topology of the complement of a ball). Therefore, under what suitable assumptions that the BCN Conjecture is true is investigated subsequently. The rigidity result in unbounded $\Omega$ is difficult because neither the moving plane method nor the method employing some integral identities can be applied directly. The reason is that the comparison principle cannot be used in unbounded domains. Moreover, the integral identities such as Pohozaev's identity cannot hold for unbounded domains.

In 1997, Reichel \cite{MR1463801} proved the existence of a bounded solution in an exterior domain $\Omega$ if $f$ satisfies some certain conditions implies $\Omega$ is the complement of a ball. For the BCN Conjecture in dimension $2$, Ros and Sicbladi \cite{MR3062759} showed it is true if the domain $\Omega$ is narrow or $f(u)\geq \lambda u$ for some $\lambda>0$. Moreover, Ros, Ruiz and Sicbaldi \cite{MR3666566} proved it holds if $\Omega$ is diffeomorphic to a half-plane and $\kappa\neq 0$. Rigidity results in unbounded $\Omega$ in other dimensions were obtained by Farina and Valdinoci \cite{MR2680184, MR2591980} and Wang and Wei \cite{MR3952780} under the assumptions that $\Omega$ is a globally Lipschitz epigraph. In 2024, Ruiz, Sicbaldi and Wu \cite{Ruiz_Sicbaldi_Wu24} presented a rigidity result in $\mathbb{R}^n$ that if the primitive $F\in C^2(\mathbb{R})$ of $f$ and $\kappa$ satisfies $\kappa^2+2F(0)\geq 0$, then either the mean curvature of $\partial \Omega$ is negative or $\Omega$ is a half-space.


\medskip

In this paper, we establish a rigidity result for overdetermined $p$-Laplace problem \cref{e.P}. Specifically, we prove that under a proper structural condition on $f$, the domain $\Omega$ must satisfy a geometric dichotomy: either the mean curvature $H(q)$ of $\partial \Omega$ is strictly negative for any boundary point $q$ or $\Omega$ is necessarily a half-space. This geometric conclusion is formalized in the following theorem:
\begin{theorem}\label{coro}
Let $\Omega \subset\mathbb{R}^{n}$ be a $C^1$ domain, $f$ be a $C^1(\mathbb{R})$ function and $\kappa \neq 0$. Let $u$ be a bounded $C^3$ solution of \eqref{e.P}. Suppose that there exists a nonpositive primitive $F\in C^{2}(\mathbb{R})$ of $f$ such that
\begin{equation}\label{cond1}
 F(0) \geq -\frac{p-1}{p}\kappa^p.
\end{equation}
For $p>2$, we assume additionally that if $F(u_0)=0$,
\begin{equation}\label{e1.1}
F(u)=O(|u-u_0|^p)~~\mbox{ as }~~u\rightarrow u_0.
\end{equation}
Then either $\Omega$ is a half-space and $u$ is parallel, or $H(q) < 0$ for any $q \in \partial\Omega$, where $H(q)$ is the mean curvature of $\partial \Omega$ at $q$.
\end{theorem}

\begin{remark}\label{Re.1}
We say that a solution $u$ is parallel (or $1$-dimensional) if there exist $x_0, a \in \mathbb{R}^n$ and a function $g: [0, +\infty) \rightarrow\mathbb{R}$ such that
\begin{equation*}
  u(x)=g\left(a\cdot (x-x_0)\right),~~~~~ x \in \{x \in \R^n:  \ a \cdot (x-x_0) >0 \}.
\end{equation*}
\end{remark}

If $\Omega$ is a bounded domain, there must exist one point on $\partial \Omega$ where the mean curvature is nonnegative. Based on this fact, we have the following corollary.
\begin{corollary}\label{co1.1}
Let $\Omega,f,\kappa,u$ be as in \Cref{coro}. Suppose that $\Omega$ is a bounded domain. Then for $1<p\leq 2$, there doesn't exist a nonpositive primitive $F$ such that \eqref{cond1} holds; for $p>2$, there doesn't exist a nonpositive primitive $F$ such that both \eqref{cond1} and \eqref{e1.1} hold.
\end{corollary}

\medskip

Theorem \ref{coro} will be obtained  as a corollary of more general results, which establish Modica type estimates for \cref{e.P}. In this framework, given a primitive $F$ of $f$, we define the $Q$-function (generally called $P$-function, here $Q$-function is used in order to distinguish from the $p$ in $p$-Laplacian):
\begin{equation}\label{Pfunction}
Q(x)=\frac{p-1}{p}\cdot |\nabla u|^{p}+F(u(x)).
\end{equation}

\medskip

In 1985, Modica \cite{MR803255} provided a gradient bound estimate to semilinear elliptic Laplacian equations in whole $\mathbb{R}^n$. Specifically, if $F$ is nonpositive and $u : \mathbb{R}^{n} \rightarrow \mathbb{R}$ is a bounded $C^{3}$ solution of $\Delta u + f(u) = 0$, then $Q \leq 0$. This estimate, called Modica estimate, is dedicated to treat gradient bound in unbounded domains with the aid of $Q$-function. In 1994, Caffarelli, Garofalo and Segala \cite{MR1296785} considered more general operators including $p$-Laplacian and obtained Liouville type theorems (i.e., rigidity of solutions) in whole space based on the gradient bound estimate. Ratto and Rigoli \cite{MR1359724} extended the analogous results of Poisson equation on complete Riemannian manifolds.

\medskip

Additionally, Farina and Valdinoci \cite{MR2680184} used the Modica type estimate to prove a pointwise gradient bound for bounded solutions of $\Delta u+f(u)=0$ in epigraphs with boundary of nonnegative mean curvature and this result is also extended to compact manifolds with nonnegative Ricci tensor, see \cite{MR2812957, MR4211208}. For more general unbounded domains with unbounded boundary, Modica type estimates of $\Delta u+f(u)=0$ have been studied in \cite{Ruiz_Sicbaldi_Wu24}. Inspired by this result, we propose to understand if these results can also hold for the problem \cref{e.P}. In other words, we aim to get the following result.

\medskip

\begin{theorem} \label{Th.P1}
Let $\Omega \subset\mathbb{R}^{n}$ be a $C^1$ domain, $f$ be a $C^1(\mathbb{R})$ function and $u$ be a bounded $C^3$ solution of \eqref{e.P}. Let $Q$ be given by \cref{Pfunction}, where $F\in C^{2}(\mathbb{R})$ is a nonpositive primitive of $f$. Then
\begin{equation}\label{Modest}
  Q(x)\leq\max\{0,F(0)+\frac{p-1}{p}\kappa^p\} \ \mbox{ for all } x \in \Omega.
\end{equation}

Moreover, if there exists a point $x_{0}\in\Omega$ such that
\begin{equation*}
  Q(x_{0})=\max\{0,F(0)+\frac{p-1}{p}\kappa^p\},
\end{equation*}
(we also assume \cref{e1.1} for $p>2$) then $Q$ is constant, $\Omega$ is a half-space and $u$ is parallel.
\end{theorem}

\begin{remark}\label{re2.1}
In this theorem, the constant $\kappa$ can be $0$.
\end{remark}

\begin{remark}\label{re2.2}
Since $u$ is a bounded positive solution, there exists $L>0$ such that $0\leq u\leq L$. Note that we assume that $f,F$ are smooth in the whole $\mathbb{R}$ (i.e., $f\in C^1(\mathbb{R}),F\in C^2(\mathbb{R})$) rather than being smooth in the range of $u$ (i.e., $f\in C^1([0,L]),F\in C^2([0,L])$). The reason is that we need to use the fact that if $F$ attains its maximum at some point $u_0$, then $F'(u_0)=0$.
\end{remark}

\medskip

When $Q$ is bounded above by $F(0)+(p-1)\kappa^p/ p,$ we can derive results concerning the mean curvature of $\partial \Omega$, as demonstrated in the following statement:

\begin{theorem}\label{Th.P2}
Let $\Omega \subset\mathbb{R}^{n}$ be a $C^1$ domain, $f$ be a $C^1(\mathbb{R})$ function and $u$ be a bounded $C^3$ solution of \eqref{e.P}. Assume that $\kappa \neq 0$ and
\begin{equation} \label{cond}
Q(x) \leq F(0)+\frac{p-1}{p}\kappa^p \ \mbox{ for all } x \in \Omega.
\end{equation}
Then $H(q) \leq 0$ for any $q \in \partial\Omega$. Moreover, if there exists $q \in \partial \Omega$ such that $H(q)=0$, then $Q$ is constant, $\Omega$ is either a half-space or a slab and $u$ is parallel.
\end{theorem}

\begin{remark}\label{Re.2}
We note that a slab is the part of the space between two parallel hyperplanes, i.e. $\mathcal{S} = \{x \in \R^n:  \ a \cdot (x-x_0) \in (z_1,z_2) \}$
for some $x_0,a \in \R^n$ and $z_1, z_2 \in \R$. A parallel solution in a slab is then a solution of the form
\begin{equation*}
  u(x)=g(a\cdot (x-x_0)), \ \ x \in \mathcal{S}.
\end{equation*}
for some function $g: (z_1,z_2) \to \R$.
\end{remark}

\section{Proof of Main Results}\label{sec2}
In this section, we first prove \Cref{Th.P1} and then we will deduce \Cref{Th.P2} and from that \Cref{coro} follows. Suppose that the assumptions of \Cref{Th.P1} are satisfied. We start with an important lemma which shows that the $Q$-function is a subsolution of some elliptic PDE. It allows us to use the maximum principle. This fact has been proved in \cite[Theorem 2.2]{MR1296785} for a general equation with divergence forms and we present here another different proof for the sake of completeness.

In the following we will write $f$ instead of $f(u)$, as well as $f'$ instead of $f'(u)$. We use $x=(x_1,...,x_n)$ to denote a general vector in $\mathbb{R}^n$. In addition, $u_i$ denotes $\partial u/\partial x_i$ and $u_{ij}$ denotes $\partial^2u/\partial x_i\partial x_j$ etc. ($1\leq i,j\leq n$). We also use the usual Einstein summation convention dropping the symbol of sum, that is understood when indices are repeated.

\begin{lemma}\label{Subsolution}
The function $Q$ satisfies:
\begin{equation*}
 \Delta Q+\frac{(p-2) u_iu_j}{|\nabla u|^{2}} Q_{ij}
  +\frac{p}{|\nabla u|^p}fu_iQ_i \geq 0,
\end{equation*}
at any $x\in \Omega$ with $\nabla u(x)\neq 0$.
\end{lemma}

\begin{proof}
Note that $Q\in C^2$ since $u\in C^3$. By differentiating $Q$, we have
\begin{equation}\label{P}
    Q_i=(p-1)|\nabla u|^{p-2}u_k u_{ki}+fu_i
\end{equation}
and
\begin{equation}\label{P-2}
  \begin{aligned}
Q_{ij}=&(p-1)|\nabla u|^{p-2}\left(u_k u_{ijk}+u_{ki}u_{kj}\right)
    +(p-1)(p-2)|\nabla u|^{p-4}u_{k}u_{ki}u_{l}u_{lj}\\
    &+f'u_iu_j+f u_{ij}.
  \end{aligned}
\end{equation}
Obviously,
\begin{equation}\label{P-3}
  \begin{aligned}
  \Delta Q=&(p-1)|\nabla u|^{p-2}\left(u_k \Delta u_{k}+u_{ki}u_{ki}\right)
    +(p-1)(p-2)|\nabla u|^{p-4}u_{k}u_{ki}u_{l}u_{li}\\
    &+f'|\nabla u|^2+f\Delta u.
  \end{aligned}
\end{equation}
In addition, we have
\begin{equation*}\label{P-form-1}
Q_iu_i=(p-1)|\nabla u|^{p-2}u_iu_k u_{i k}+f|\nabla u|^2.
\end{equation*}
Then
\begin{equation}\label{p-form-2}
  u_iu_k u_{i k}=\frac{Q_iu_i-f|\nabla u|^2}{(p-1)|\nabla u|^{p-2}}
\end{equation}
holds at the point $x$ that $\nabla u(x) \neq 0$.

Next, by \cref{P} and applying the Cauchy-Schwarz inequality to $u_k u_{ki}$, we get
\begin{equation*}
\begin{aligned}
|Q_i-fu_i|=(p-1)|\nabla u|^{p-2}|u_ku_{ki}|
\leq (p-1)|\nabla u|^{p-1} \left(\sum_{k=1}^{n} u_{ki} u_{ki}\right)^{\frac{1}{2}}.
\end{aligned}
\end{equation*}
Then
\begin{equation}\label{estimate-uki}
\begin{aligned}
|\nabla u|^{p-2} u_{ki}u_{ki}\geq & \frac{1}{(p-1)^2|\nabla u|^p}\left(|\nabla Q|^2-2fu_iQ_i+f^2|\nabla u|^2\right)
\end{aligned}
\end{equation}
for any $x$ that $\nabla u(x) \neq 0$.
Write now the equation in \cref{e.P} in non-divergence form:
\begin{equation}\label{eq-div}
  |\nabla u|^{p-2}\Delta u+(p-2)|\nabla u|^{p-4}u_iu_ju_{ij}+f(u)=0.
\end{equation}
Then, by differentiating \cref{eq-div} with respect to $x_k$, we get
\begin{equation*}
\begin{aligned}
  |\nabla u|^{p-2}\Delta u_k+(p-2)|\nabla u|^{p-4} u_{l}u_{lk}\Delta u
  +(p-2)|\nabla u|^{p-4}\left(2u_{ik}u_ju_{ij}+u_iu_ju_{ijk}\right)\\
  +(p-2)(p-4)|\nabla u|^{p-6}u_iu_ju_{ij}u_{l}u_{lk}+f' u_k=0.
\end{aligned}
\end{equation*}
Next, multiply by $u_k$ both sides in above equation and take the sum over $k$:
\begin{equation}\label{diff-eq}
\begin{aligned}
  |\nabla u|^{p-2}u_k\Delta u_k
  +(p-2)|\nabla u|^{p-4}\left(u_k u_{l}u_{lk}\Delta u+2u_ku_{ik}u_ju_{ij}+u_iu_ju_ku_{ijk}\right)\\
  +(p-2)(p-4)|\nabla u|^{p-6}u_iu_ju_{ij}u_ku_{l}u_{lk}+f' |\nabla u|^2=0.
\end{aligned}
\end{equation}
Inserting \cref{p-form-2} into \cref{eq-div}, we obtain
\begin{equation}\label{Laplace-u}
\Delta u=-\frac{p-2}{(p-1)|\nabla u|^{p}}Q_i u_i-\frac{1}{(p-1)|\nabla u|^{p-2}} f.
\end{equation}
By considering \cref{P-3}$/(p-1)+$\cref{P-2}$\times (p-2) u_iu_j/\left((p-1)|\nabla u|^{2}\right)-$\cref{diff-eq}, we have
\begin{equation}\label{P-eq}
\begin{aligned}
  &\frac{1}{p-1}\Delta Q+\frac{(p-2) u_iu_j}{(p-1)|\nabla u|^{2}} Q_{ij}\\
=&|\nabla u|^{p-2} u_{ki}^2+\frac{f\Delta u}{p-1}+\frac{(p-2)f u_iu_ju_{ij}}{(p-1)|\nabla u|^{2}}\\
&+2(p-2)|\nabla u|^{p-6}\left(u_iu_ju_{ij}\right)^2
-(p-2)|\nabla u|^{p-4}u_i u_{j}u_{ij}\Delta u.
\end{aligned}
\end{equation}
There are five terms in the righthand. For the second and fifth terms, use \cref{Laplace-u} to replace $\Delta u$ and for the third to fifth terms, use \cref{p-form-2} to replace $u_i u_j u_{ij}$. Then we have
\begin{equation}\label{P-eq4}
  \begin{aligned}
    \frac{f\Delta u}{p-1}=&-\frac{(p-2)f Q_iu_i}{(p-1)^2|\nabla u|^p}-\frac{f^2}{(p-1)^2|\nabla u|^{p-2}};\\
    \frac{(p-2)f u_iu_ju_{ij}}{(p-1)|\nabla u|^{2}}=&\frac{(p-2)f Q_iu_i}{(p-1)^2|\nabla u|^p}-\frac{(p-2)f^2}{(p-1)^2|\nabla u|^{p-2}};\\
    2(p-2)|\nabla u|^{p-6}\left(u_iu_ju_{ij}\right)^2=&\frac{2(p-2)Q_iu_iQ_ju_j}{(p-1)^2|\nabla u|^{p+2}}
    -\frac{4(p-2)f Q_iu_i}{(p-1)^2|\nabla u|^p}+\frac{2(p-2)f^2}{(p-1)^2|\nabla u|^{p-2}};\\
    -(p-2)|\nabla u|^{p-4}u_i u_{j}u_{ij}\Delta u=&\frac{(p-2)^2 Q_iu_iQ_ju_j}{(p-1)^2|\nabla u|^{p+2}}
    -\frac{(p-2)(p-3)f Q_iu_i}{(p-1)^2|\nabla u|^p}-\frac{(p-2)f^2}{(p-1)^2|\nabla u|^{p-2}}.
  \end{aligned}
\end{equation}
Next, by inserting \cref{P-eq4} and \cref{estimate-uki} into \cref{P-eq}, we obtain
\begin{equation*}
\begin{aligned}
  &\frac{1}{p-1}\Delta Q+\frac{(p-2) u_iu_j}{(p-1)|\nabla u|^{2}} Q_{ij}\\
=&|\nabla u|^{p-2} u_{ki}^2+\frac{p(p-2)Q_iu_iQ_ju_j}{(p-1)^2|\nabla u|^{p+2}}
-\frac{(p-2)(p+1)fQ_i u_i}{(p-1)^2|\nabla u|^p} -\frac{f^2}{(p-1)^2|\nabla u|^{p-2}}\\
\geq & \frac{|\nabla Q|^2}{(p-1)^2|\nabla u|^p}
+\frac{p(p-2)Q_iu_iQ_ju_j}{(p-1)^2|\nabla u|^{p+2}}
-\frac{pfQ_iu_i}{(p-1)|\nabla u|^p}.
\end{aligned}
\end{equation*}
That is
\begin{equation*}
  \frac{1}{p-1}\Delta Q+\frac{(p-2) u_iu_j}{(p-1)|\nabla u|^{2}} Q_{ij}
  +\frac{pfu_iQ_i}{(p-1)|\nabla u|^p}
\geq \frac{|\nabla Q|^2}{(p-1)^2|\nabla u|^p}
+\frac{p(p-2)Q_iu_iQ_ju_j}{(p-1)^2|\nabla u|^{p+2}} .
\end{equation*}
If $p\geq 2$, the last term is nonnegtive and the conclusion follows. If $1<p<2$, by the Cauchy-Schwarz inequality,
\begin{equation*}
  |\nabla Q|^2+p(p-2)\frac{Q_iu_iQ_ju_j}{|\nabla u|^2}
  \geq |\nabla Q|^2(1+p(p-2))=|\nabla Q|^2(p-1)^2\geq 0.
\end{equation*}
This concludes the proof of the lemma.
\end{proof}

Next, we provide with a uniform estimate of the gradient of $u$, which will be essential in the proof of \Cref {Th.P1} and thus \Cref{coro}. Throughout the paper, we set $B_r(x_0):=\{x\in \mathbb{R}^n:|x-x_0|<r \}$ and $B_r=B_r(0)$ for simplicity. In addition, we denote $\Delta_p u:=\mathrm{div} \left(|\nabla u|^{p-2}\nabla u\right)$ in the following.

\begin{lemma} \label{Le21} If $u$ is a bounded solution of the problem \Cref{e.P}, then there exists a constant $M>0$ depending only on $p,n$, $\|u\|_{L^{\infty}}$,  $\| f(u) \|_{L^{\infty}}$ and $\kappa$ such that
\begin{equation*}
   |\nabla u|\leq M~~\mbox{ in }~~\Omega.
\end{equation*}
\end{lemma}
\begin{proof}
Define
\begin{equation*}
  \Omega':=\{x\in\Omega:\mbox{dist}(x,\partial\Omega)>1\}.
\end{equation*}
By the interior regularity (see \cite[Theorem 1]{MR709038} or \cite[Theorem 1]{MR727034}),
\begin{equation}\label{e2.3}
   |\nabla u|\leq M'~~\mbox{ in }~~\Omega',
\end{equation}
where $M'$ depends only on $p$, $n$, $\|u\|_{L^{\infty}}$ and $\| f(u) \|_{L^{\infty}}$.

\medskip

In the following, we prove the gradient bound in $\Omega\setminus \Omega'$. Given any $x_{0} \in \Omega\setminus \Omega'$, denote by $h:=\mbox{dist}(x_{0},\partial\Omega)\leq 1$ and assume that this distance is attained at $y_{0}\in\partial\Omega$. Then we set
\[\tilde{u}(x):=\frac{1}{h}u(x_{0}+hx).\]
This gives
\begin{align*}
  \Delta_p\tilde{u} & =|\nabla \tilde{u}|^{p-2}\Delta \tilde{u}+(p-2)|\nabla \tilde{u}|^{p-4}\tilde{u}_i\tilde{u}_j\tilde{u}_{ij} \\
  &=h|\nabla u|^{p-2}\Delta u+h(p-2)|\nabla u|^{p-4}u_iu_ju_{ij}\\
  &= h\Delta_p u.
\end{align*}
Then we have
\begin{equation*}
  \Delta_p\tilde{u}+hf(h\tilde{u})=0~~\mbox{ and }~~\tilde{u}>0~~\mbox{ in }~~B_1.
\end{equation*}

\medskip

We claim that
\begin{equation*}
  \sup\limits_{B_{1/2}}\tilde{u}(x)\leq M,
\end{equation*}
where $M$ is a large constant to be specified later.

Let $z_{0}:=\left(y_{0}-x_{0}\right)/h$ and then
\begin{equation}\label{e2.1}
  |\nabla\tilde{u}(z_{0})|=|\nabla u(y_{0})|=\kappa.
\end{equation}
By applying the Harnack inequality (see \cite[Theorem 5, Theorem 6 and Theorem 9]{MR170096}), one can get that
\begin{equation}\label{harnack}
\sup\limits_{B_{1/2}} \tilde{u}(x)\leq C_0\left(\inf\limits_{B_{1/2}}
\tilde{u}(x)+h\|f(u)\|_{L^{\infty}(B_{1})}\right),
\end{equation}
where $ C_0=C_0(n,p).$

Reasoning by contradiction, we suppose that
\[\sup\limits_{B_{1/2}}\tilde{u}(x)> M.\]
It follows from \cref{harnack} and taking $M$ large enough that
\begin{equation*}
  \inf\limits_{B_{1/2}} \tilde{u}(x)
  >\frac{M}{C_0}-h\|f(u)\|_{L^{\infty}(B_{1})}>\frac{M}{2C_0}.
\end{equation*}

We now denote a function
\begin{equation*}
  v(x)=C_1(|x|^{-\beta}-1)
\end{equation*}
with
\begin{equation*}
  C_1= \frac{M}{2C_0(2^{\beta}-1)}.
\end{equation*}
By choosing $\beta$ large enough first (depending only on $p, n$) and then choosing $M$ large enough (depending on $\beta$), we have
\begin{equation*}
 \begin{cases}
 \Delta_p v\geq \|f(u)\|_{L^{\infty}} &\mbox{in $B_1\setminus B_{1/2}$},\\
v\leq \tilde{u}&\mbox{on $\partial B_{1/2}$},\\
 v=0&\mbox{on $ \partial B_{1}$}.
\end{cases}
\end{equation*}

Then by comparison principle \cite[Theorem 1.2]{MR1632933}, we have
\[\tilde{u}\geq v \qquad \mbox{in \,\, $B_1\setminus B_{1/2}$}.\]
By a direct computation,
 \[v\geq C C_1\mbox{dist$(x,\partial\Omega)$}\geq CM\mbox{dist$(x,\partial\Omega)$}.\]
 Then
 \[\tilde{u}\geq CM \mbox{dist$(x,\partial\Omega)$},\]
i.e.,
 \[|\tilde{u}_\nu(z_0)|\geq CM.\]
 By noting \cref{e2.1},
 \begin{equation*}
\kappa\geq CM,
 \end{equation*}
which is a contradiction if we take $M$ large enough. Therefore,
\[\sup\limits_{B_{1/2}}\tilde{u}(x)\leq M.\]
By the interior gradient estimate again,
\begin{equation}\label{e2.2}
|\nabla u(x_{0})|=|\nabla\tilde{u}(0)|\leq C M,
\end{equation}
where $C$ depends only on $p,n$ and $\|f(u) \|_{L^{\infty}}$.

Finally, from \cref{e2.3} and \cref{e2.2}, we arrive at the conclusion.
\end{proof}

\begin{remark}\label{re.2.2}
In \cite[Lemma 2.2]{Ruiz_Sicbaldi_Wu24}, the authors decompose $\tilde{u}$ into the sum of two solutions of two equations. This method cannot be used in nonlinear equations directly. We use a barrier function here to obtain the boundedness which is quite flexible. Even for $p=2$, this proof for gradient bound is simpler.
\end{remark}

Now, we are ready to present the proof of our main result Theorem \ref{Th.P1}. The proof will be shown by two propositions. Define
\begin{equation}\label{alpha}
  \hat \alpha :=\max\{0,F(0)+\frac{p-1}{p}\kappa^{p}\}.
\end{equation}
First, we prove the inequality \cref{Modest}, which is based on \cite[Theorem 1.6]{MR1296785} and \cite[Proposition 3.1]{Ruiz_Sicbaldi_Wu24}.

\begin{proposition}\label{pr2.1}
Under the assumptions of Theorem \ref{Th.P1}, we have that $Q(x) \leq \hat \alpha$ for all $x \in \Omega$. \end{proposition}

\begin{proof}
According to Lemma \ref{Le21}, $Q$ is bounded. We prove the proposition by contradiction. Suppose that
\[\beta:=\sup\limits_{\Omega} Q(x)>\hat \alpha.\]
Then there exists a sequence $\{x_{m}\}_{m\in\mathbb{N}}\subset\Omega$ such that
\[Q(x_{m})\rightarrow\beta\,\ \mbox{as}\,\ m\rightarrow \infty.\]
If $x_{m}$ is bounded, up to a subsequence, we can assume that $x_m \to x_{0}$ with $Q(x_0)=\beta$. Notice that $\nabla u(x_0) \neq 0$ since $\beta>0$ and $F\leq 0$. Moreover, $x_0 \in \Omega$ since $\beta > F(0)+(p-1)\kappa^p /p $. As a consequence, $Q$ attains a local maximum at an inner point $x_0$. Then, by the strong maximum principle, $Q$ is constant on a neighborhood of $x_0$. This argument implies that the set:
 $$ \{x \in \Omega: \ Q(x)= \beta\}$$
is a non-empty open subset of $\Omega$. Since it is obviously closed, then $Q(x)\equiv \beta$ for all $x \in \Omega$ which contradicts the fact that $Q(x)= F(0)+(p-1)\kappa^p /p < \beta$ if $x \in \partial \Omega$.

\medskip

If $x_{m}$ is unbounded, we will discuss it in two cases.

\medskip

\textbf{Case 1:} $\limsup_{m \to +\infty} \mbox{dist}(x_{m},\partial\Omega)>\delta$ for some positive constant $\delta$.

Let us extend $u$ by $0$ outside $\Omega$, so that $u: \R^n \to \R$ is a globally Lipschitz function. Consider $u_{m}(x):=u(x+x_{m})$ and $Q_{m}(x):=Q(x+x_{m})$, then
\begin{equation*}
Q_{m}(0)=Q(x_{m})\to\beta.
\end{equation*}
Take a subsequence, still denoted by $u_m$, so that on compact sets of $\R^n$ the sequence $u_m$ converges uniformly to a certain Lipschitz function $u_\infty \geq 0$. Let
\begin{equation*}
  \Omega_{\infty}= \{x \in \R^n: \ u_\infty (x) >0\}
\end{equation*}
and
\begin{equation*}
  Q_{\infty}= \frac{p-1}{p}|\nabla u_{\infty}|^p +  F(u_\infty) \in L^{\infty}(\R^n).
\end{equation*}
We first claim that $\Omega_{\infty}$ is not empty. Observe that $\mbox{dist}(x_{m},\partial\Omega)>\delta$, we have
\begin{equation}\label{um}
  \Delta_p u_m+f(u_m)=0~~\mbox{ in }~~B_{\delta}.
\end{equation}
Since $f(u_m)$ is a globally Lipschitz function in $\R^n$ and $\nabla u_m$ is uniformly bounded (according to \Cref{Le21}), by interior $C^{1,\alpha}$ regularity (see \cite{MR709038, MR721568, MR727034}), we have that $u_m$ is bounded in $C^{1, \alpha}$ on any compact subset of $B_{\delta}$. By taking the limit in \cref{um}, $u_{\infty} \in C^{1,\alpha}$ is a weak solution of
\begin{equation}\label{u-infty}
  \Delta_p u_{\infty}+f(u_{\infty})=0~~\mbox{ in }~~B_{\delta}.
\end{equation}
Moreover, note that $F\leq 0$ and
\begin{equation*}
Q_{\infty}(0)=\lim_{m\to \infty} Q_m(0)=\beta>0.
\end{equation*}
Then $\nabla u_{\infty}(0)\neq 0$. This excludes the possibility $u_{\infty}(0)=0$. Indeed, if $u_{\infty}(0)=0$ then $u_{\infty}$ attains the minimal point at $0$ and $\nabla u_{\infty}(0)$ will be $0$. Thus, $\Omega_{\infty}$ is not empty since $0 \in \Omega_{\infty}$.

For any $\Omega''\subset \subset \Omega'\subset\subset \Omega_{\infty}$, there exists $a>0$ such that $u_{\infty}\geq a$ in $\Omega'$. Because $u_m\rightarrow u_{\infty}$ uniformly, one can conclude that $u_m>0$ in $\Omega'$ for $m$ large enough. Thus, $u_m$ satisfies \cref{um} in $\Omega'$. Similarly, by the interior $C^{1,\alpha}$ regularity, we have that $u_m\rightarrow u_{\infty}$ in $C^{1,\alpha}(\Omega'')$. Hence, $u_{\infty}$ satisfies \cref{u-infty} in $\Omega_{\infty}$ and $Q_m\rightarrow Q_{\infty}$ in any compact set of $\Omega_{\infty}$. Then we get
\begin{equation*}
  Q_{\infty} \leq \beta~~\mbox{ in }~~\Omega_{\infty}.
\end{equation*}
Recall that $Q_{\infty}(0)=\beta$. By applying the strong maximum principle to $Q_{\infty}$, we conclude that $Q_{\infty}$ is constant on a neighborhood of $0$. Thus, we have that the set $\{ x \in \Omega:\ Q_{\infty}(x)= \beta \}$ is an open subset of $\Omega_{\infty}$. In addition, by the continuity of $Q_{\infty}$ in $\Omega_{\infty}$, it is also closed. As a consequence,
\begin{equation*}
 Q_{\infty}(x) = \beta \mbox{ for all } x \in \tilde{\Omega}_{\infty}, \end{equation*} where $\tilde{\Omega}_{\infty}$ is the connected component of $\Omega_{\infty}$ containing $0$.

Since $u_{\infty}$ is bounded, then there exists a sequence $y_m\in \tilde{\Omega}_{\infty}$ such that $\nabla u_{\infty}(y_m)\rightarrow 0$. Then
\begin{equation*}
  0<\beta=\lim_{m\to\infty}Q_{\infty}(y_m)
  \leq \frac{p-1}{p}\lim_{m\to\infty}|\nabla u_{\infty}(y_m)|^p
  +\limsup_{m\to\infty}F(u_{\infty}(y_m))\leq 0,
\end{equation*}
which is a contradiction.

\medskip

\textbf{Case 2:}  $\lim_{m \to +\infty} \mbox{dist}(x_{m},\partial\Omega)=0.$

As in case 1, we consider $u$ extended by $0$ outside $\Omega$. Then $u: \R^n \to \R$ is a globally Lipschitz function. Observe that $u(x_{m})\rightarrow 0$ and $Q(x_m)\rightarrow \beta$. Since $|\nabla u|$ is bounded by Lemma \ref{Le21}, then we have
\begin{equation} \label{jo3}
\frac{p-1}{p}|\nabla u(x_{m})|^{p} = Q(x_m) -  F(u(x_m)) \to \beta -F(0) >\frac{p-1}{p}\kappa^{p}.
\end{equation}
Denote $h_{m}:=\mbox{dist}(x_{m},\partial\Omega)\rightarrow0$ and assume that this distance is attained at $y_{m}\in\partial\Omega.$ Let
\[u_{m}(x):=\frac{1}{h_{m}}u(y_{m}+h_{m}x).\]
Then, $u_m(0)=0$ and $\nabla u_m=\nabla u$ is uniformly bounded. Therefore, we can take a subsequence such that $u_m$ converges uniformly in compact sets to a limit function $u_\infty \geq 0$, which is Lipschitz continuous (but possibly unbounded).

Let us consider $\Omega_{\infty}= \{x \in \R^n: \ u_{\infty}(x) >0\}$ as before. First, we show that $\Omega_{\infty}$ is non-empty. Let $z_{m}:=\left(x_{m}-y_{m}\right)/h_{m}$ and we can assume that $z_{m}$ converges to $z_{\infty}$ since $|z_{m}|=1$. Thus, since $\nabla u_m$ is uniformly bounded and $u_m$ satisfies
\begin{equation*}
  \Delta_p u_m +h_m f(u_m)=0~~\mbox{ in }~~B_1(z_m).
\end{equation*}
By the interior $C^{1,\alpha}$ regularity, we obtain that $u_m\rightarrow u_{\infty}$ in $C^{1,\alpha}$ sense in compact sets of $B_1(z_{\infty})$ and $u_{\infty}$ satisfies
\begin{equation}\label{u-infty-case2}
\Delta_p u_{\infty} =0~~\mbox{ in }~~ B_1(z_{\infty}).
\end{equation}
Note that $u_{\infty}\geq 0$ in $B_1(z_{\infty})$ and according to \cref{jo3},
\begin{equation}\label{Du-infty-case2}
a:=|\nabla u_{\infty}(z_{\infty})|=\lim_{m\rightarrow \infty} |\nabla u_{m}(z_{m})|
=\lim_{m\rightarrow \infty} |\nabla u(x_{m})|>\kappa.
\end{equation}
By the strong maximum principle we can obtain $u_{\infty}>0$ in $B_1(z_{\infty})$. It means that $B_1(z_{\infty}) \subset \Omega_{\infty}$ and $\Omega_{\infty}$ is non-empty. Similarly, the same argument works in $\Omega_{\infty}$ and then we have that $u_m \rightarrow u_{\infty}$ in $C^{1,\alpha}$ sense in compact sets of $\Omega_{\infty}$. Now $u_{\infty}$ satisfies \cref{u-infty-case2} in $\Omega_{\infty}$.

Next, we claim that $|\nabla u_{\infty}|^p$ is a subsolution of an elliptic PDE. First, let's rewrite \cref{u-infty-case2} as follows.
\begin{equation}\label{u-inf2}
  |\nabla u_{\infty}|^{p-2}\left(\delta^{ij}+(p-2)\frac{(u_{\infty})_i (u_{\infty})_j}{|\nabla u_{\infty}|^2}\right)
  (u_{\infty})_{ij}=0.
\end{equation}
Let
\begin{equation*}
  \begin{aligned}
    a^{ij}&=\delta^{ij}+(p-2)\frac{(u_{\infty})_i (u_{\infty})_j}{|\nabla u_{\infty}|^2},\\
    A^{ij}&=|\nabla u_{\infty}|^{p-2} a^{ij}.
  \end{aligned}
\end{equation*}
Then, \cref{u-inf2} becomes
\begin{equation*}
  A^{ij}(u_\infty)_{ij}=0.
\end{equation*}
Next, by taking derivatives of \cref{u-infty-case2} with respect to $x_k$, we have
\begin{equation*}
  \left(|\nabla u_{\infty}|^{p-2}(u_{\infty})_{ik}
    +(p-2)|\nabla u_{\infty}|^{p-4}(u_{\infty})_{i}(u_{\infty})_{j}(u_{\infty})_{jk}\right)_i=0,
\end{equation*}
which can be rewritten as
\begin{equation*}
  \left(|\nabla u_{\infty}|^{p-2}\left(\delta^{ij}+(p-2)\frac{(u_{\infty})_i (u_{\infty})_j}{|\nabla u_{\infty}|^2}\right)(u_{\infty})_{jk}\right)_i=0.
\end{equation*}
That is,
\begin{equation}\label{u-inf-k}
  \left(A^{ij}(u_{\infty})_{jk}\right)_i=0.
\end{equation}
In the following, by a direct computation and combining \cref{u-inf-k}, we obtain
\begin{align}
   \left(a^{ij}\left(|\nabla u_{\infty}|^p\right)_j\right)_i&=\left(a^{ij}p|\nabla u_{\infty}|^{p-2}(u_{\infty})_k(u_{\infty})_{kj}\right)_i \nonumber\\
    &=p\left(\left(A^{ij}(u_{\infty})_{jk}\right)_i (u_{\infty})_k+A^{ij}(u_{\infty})_{kj}(u_{\infty})_{ki}\right)\nonumber\\
    &=pA^{ij}(u_{\infty})_{kj}(u_{\infty})_{ki}\geq 0 \label{Dup}.
  \end{align}
That is, $|\nabla u_{\infty}|^{p}$ is a subsolution of some elliptic PDE.

\medskip

Note that
\begin{equation*}
  \begin{aligned}
    \frac{p-1}{p}|\nabla u_{\infty}(z_{\infty})|^{p}+F(0)
    &=\lim_{m \rightarrow \infty} \Big(\frac{p-1}{p}|\nabla u_{m}(z_{m})|^{p}+F(h_m u_{m}(z_{m}))\Big)\\
    &=\lim_{m \rightarrow \infty} \Big(\frac{p-1}{p}|\nabla u(x_{m})|^{p}+F(u(x_{m}))\Big)\\
    &=\lim_{m \rightarrow \infty} Q(x_m)=\beta.
  \end{aligned}
\end{equation*}
Additionally, for all $ x \in B_1(z_\infty)$,
\begin{equation*}
  \begin{aligned}
    \frac{p-1}{p}|\nabla u_{\infty}(x)|^{p}+F(0)
    &=\lim_{m \rightarrow +\infty} \Big(\frac{p-1}{p}|\nabla u_{m}(x)|^{p}+F(h_m u_{m}(x))\Big)\\
    &=\lim_{m \rightarrow +\infty} \Big(\frac{p-1}{p}|\nabla u(y_m+h_m x)|^{p}+F(u(y_m+h_m x))\Big)\\
    &=\lim_{m \rightarrow +\infty} Q(y_m+h_m x)\leq \beta.
  \end{aligned}
\end{equation*}
As a consequence,
\begin{equation*}
  |\nabla u_{\infty}(z_\infty)| \geq |\nabla u_{\infty} (x)|
\end{equation*}
for all $ x \in B_1(z_\infty)$.

Since $|\nabla u_{\infty}|^p$ (see \cref{Dup}) is a subsolution, then by the strong maximum principle we obtain
\begin{equation*}
  |\nabla u_{\infty}|\equiv |\nabla u_{\infty}(z_{\infty})|=a ~~\mbox{ in }~~\tilde{\Omega}_{\infty},
\end{equation*}
where $a$ is given in \eqref{Du-infty-case2} and $\tilde{\Omega}_{\infty}$ is the connected component of $\Omega$ containing $B_1(z_{\infty})$. Thus, $u_{\infty}$ is harmonic.

Recall moreover that $u_{\infty}(0)=0$. Therefore, up to a rotation and a translation, $z_{\infty}=0$ and
\begin{equation*}
  u_{\infty}(x)=ax_{n}~~\mbox{ in }~~\tilde{\Omega}_{\infty}.
\end{equation*}
In particular, $\tilde{\Omega}_{\infty}$ is the upper half-space.

\medskip

Next, we will arrive at the contradiction. For any $\varepsilon\in(0,1)$, by the uniform convergence of $u_{m}$ to $u_{\infty}$,
\begin{equation*}
  |u_{m}-ax_{n}|<\varepsilon^2~~\mbox{ in }~~B_1^+
\end{equation*}
for all large enough $m$, where $B_1^+$ is the upper half ball and
\begin{equation*}
  \Delta_p u_m \leq \varepsilon~~\mbox{ in }~~\{x_n>\varepsilon\}\cap B_1^+.
\end{equation*}

Denote $x=(x', x_n)$ where $x'=(x_{1},\cdot\cdot\cdot,x_{n-1})$. Inspired by \cite[Lemma 4.4]{MR3567263} and \cite[Proposition 3.1]{Ruiz_Sicbaldi_Wu24}, we consider the domain $D_{t},$ the perturbation of $B^+_{1}\cap\{x_{n}>\varepsilon\}$, given by
\begin{equation*}
  D_{t}=\{x\in B^+_{1}:x_{n}>\varepsilon(\eta (x')+t),~~-1\leq t\leq 0\},
\end{equation*}
where $1\leq\eta (x')\leq 2$ is a smooth bump function defined by
\begin{equation*}
  \eta (x')=
  \begin{cases}
    1, ~~~& \mbox{in}~~ |x'|\leq\frac{1}{4}, \\
    2, ~~~& \mbox{in}~~ \frac{1}{2}\leq |x'|\leq 1.
  \end{cases}
\end{equation*}
It is clear that $u_m >0$ in  $D_{0}$ for sufficiently large $m$.

Let $\delta=(a-\kappa)/2$ and define the following barrier function
\begin{equation*}
  w(x)=(a-\delta)(x_n-\varepsilon(\eta (x')+t))+\frac{\delta}{2}(x_n-\varepsilon(\eta (x')+t))^2.
\end{equation*}
Here $w$ depends on $t$. It is worth to point out that we construct $w$ on the basis of the modification of $a x_n$. By a straightforward computation,
\begin{equation*}
  \begin{aligned}
    w_{ij}&=-\varepsilon \left(a-\delta+\delta(x_n-\varepsilon(\eta (x')+t))\right)\eta_{ij}
    +\delta \varepsilon^2 \eta_i \eta_j,~~i,j\leq n-1;\\
    w_{in}&=-\varepsilon \delta \eta_i,~~i\leq n-1;\\
    w_{nn}&=\delta.
  \end{aligned}
\end{equation*}
Observe that $\delta$ is a fixed number and by choosing $\varepsilon$ enough small, we have that $w$ satisfies
\begin{equation*}
 \begin{cases}
 \Delta_p w(x)\geq\varepsilon &\mbox{in $D_{t}$},\\
w(x)\leq u_m &\mbox{on $\partial B^+_{1}\cap\partial D_t$},\\
 w(x)=0&\mbox{on $ B^+_{1}\cap\{x_{n}=\varepsilon(\eta (x')+t)\}$},
\end{cases}
\end{equation*}
for any $-1\leq t\leq 0.$ The comparison principle implies that
\begin{equation*}
  u_m\geq w~~\mbox{ in }~~ D_t.
\end{equation*}
Since $u_m(0)=0$, there exists a $t_0\in[-1,0]$ such that $D_{t_{0}}\subseteq B_{1} \cap \{ u_m>0\}$ will touch
$$J(u_{m}):=B_{1} \cap \partial \{ u_m>0\} $$
at some point $q\in J(u_{m})\cap\{|x'|<1/2\}$. It is clear that $w(x)\leq u_{m}(x)$ in $D_{t_{0}}$. Then
\[|w_{\nu}(q)|\leq |(u_{m})_{\nu}(q)|=\kappa,\]
where $\nu$ is the outer normal to $\partial D_{t_{0}}$ at $q$. However, by the definition of $w$ and \eqref{Du-infty-case2},
\[|w_\nu(q)|\geq a-\frac{3\delta}{2}>\kappa,\]
which is a contradiction.
\end{proof}

\begin{remark}\label{re.2.3}
In the proof of \cite[Proposition 3.1]{Ruiz_Sicbaldi_Wu24}, the authors solve a function $w$ defined in $D_{t_0}$ with $\Delta w=\varepsilon$. At the touch point $p$, they estimate $|w_{\nu}(p)|$ by an approximation procedure. This method is difficult to be applied to nonlinear equations.

The approach of constructing a barrier and comparing it directly with the solution $u_m$ is more flexible. This idea is also used for the Modica type estimate for the capillary problem (see \cite[Proposition 2.5]{LianSicbaldi2025}).
\end{remark}

\begin{remark}\label{re.2.5}
This proposition is the most crucial part to obtain Modica type estimate (or the gradient bound) with regard to the unbounded domain with unbounded boundary. The main idea is to get a contradiction by showing that $Q$ attains its maximum at an interior point.
\end{remark}

It remains to show the following result to conclude the proof of Theorem \ref{Th.P1}.

\begin{proposition} \label{esta} Under the assumptions of Theorem \ref{Th.P1}, if $Q(x_0)=\hat \alpha$ at a point $x_0\in\Omega,$ then $Q$ is constant, $u$ is parallel and $\Omega$ is a half-space. \end{proposition}

\begin{proof} Proceeding as \cite[Proposition 3.2]{Ruiz_Sicbaldi_Wu24} (see also \cite{MR1296785}), we will obtain the conclusion by the following two steps.

\medskip

\textbf{Step 1:}  If $Q(x_0)=\hat \alpha$, then $\nabla u(x_0) \neq 0$.

To reach a contradiction, we assume that $\nabla u(x_{0})=0$. Then $0 \leq \hat\alpha=F(u(x_{0})) \leq 0$, so that $\hat\alpha= F(u_0)=0$ where $u_0 = u(x_0)$. Now consider the set
\[U=\{x\in \Omega:u(x)=u_{0}\}.\]
 It is obvious that $U$ is closed and $U\neq\emptyset$. Take a point $x_{1}\in U$ and consider the function $\varphi(t)=u(x_{1}+tw)-u_{0}$ for $t$ small enough, where $w\in \mathbb{S}^{n-1}$ is arbitrarily fixed. One has that
\[\varphi'(t)=\nabla u(x_{1}+tw)w.\]
Then
\begin{align*}
      |\varphi'(t)|^{p}&\leq|\nabla u(x_{1}+tw)|^{p}\leq\frac{p}{p-1}\left(\hat \alpha-F(u(x_{1}+tw))\right)\\
      &=\frac{p}{p-1}\left(F(u_{0})-F(u(x_{1}+tw))\right).
\end{align*}
By noting that $F$ attains the maximum at $u_0$, we have $F'(u_{0})=0$ and $F''(u_{0})\leq0$. Additionally, since $F(u_{0})=0$, we have $F(u)=O((u-u_{0})^{2})$ as $|u-u_{0}|\rightarrow0$. If $1<p\leq 2$,
\begin{equation*}
  |\varphi'(t)|^p\leq C|\varphi(t)|^2\leq C|\varphi(t)|^p
\end{equation*}
as $t$ small enough. If $p>2$, by the assumption \cref{e1.1}, we also obtain
\begin{equation}\label{varphi}
|\varphi'(t)|\leq C|\varphi(t)|
\end{equation}
as $t$ small enough. It follows that $\varphi\equiv0$ in $[-\varepsilon,\varepsilon]$ for some $\varepsilon>0$ since $\varphi(0)=0$. It shows that the set $U$ is open. Hence, $U=\Omega$ which implies that $u$ is a constant. This makes a contradiction since $u$ is a positive solution and $u=0$ on the boundary.

\medskip

\textbf{Step 2:}  Conclusion.

Now we have that $\nabla u(x)\neq0$ for any $x\in\Omega$ with $Q(x) =\hat \alpha$. By the maximum principle applied to $Q$, we conclude that the set
\begin{equation*}
  \{ x \in \Omega: \ Q(x)=\hat \alpha\}
\end{equation*}
is a non-empty open set, which is also closed by continuity. Then $Q(x)= \hat \alpha$ for all $x \in \Omega$. In particular, $u$ has no critical points in $\Omega$.

Let  $$v=G(u)=\left(\frac{p-1}{p}\right)^{\frac{1}{p}}\int_{u_{0}}^{u}(\hat\alpha-F(s))^{-\frac{1}{p}}ds,$$
for some fixed $u_{0}\in u(\Omega).$  First observe that the integral defining $G$ is not singular since $\hat\alpha > F(u)$. By a straightforward computation, one has
\begin{align} \label{otra}
   \Delta_p v&=|G'(u)|^{p-2}\left(G''(u)(p-1)|\nabla u(x)|^{p}+G'(u)\Delta_p u\right)\nonumber\\
   & =|G'(u)|^{p-2}\left(G''(u)p(\hat\alpha-F(u))-G'(u)F'(u)\right)\nonumber\\
   &=0.
\end{align}
Moreover,
\begin{equation}\label{eq42}
    |\nabla v|^{p}=G'(u)^{p}|\nabla u|^{p}=1.
\end{equation}
We can infer that $v(x)=a\cdot x+b$ for some $a\in \mathbb{R}^{n}$ with $|a|=1$ and $b\in \mathbb{R}$. Indeed, by \cref{otra} and \cref{eq42}, we have
\begin{equation*}
  \Delta v=\mbox{div}(|\nabla v|^{p-2}\nabla v)=\Delta_p v=0.
\end{equation*}
Fixing a point $q\in\Omega$, there exists a unit vector, say $e_1$, such that
\begin{equation*}
  v_1(q):=\nabla v \cdot e_1=1.
\end{equation*}
Thus, $v_1$ attains the maximum at $q$. Since $v_1$ is harmonic, by the strong maximum principle, $v_1\equiv 1$. By combining with $|\nabla v|\equiv 1$, we have $v_i\equiv 0$ for $2\leq i\leq n$. Therefore, $v(x)=x_{1}+b$ where $b$ is a constant.

Since $G$ is invertible, we can obtain that $u(x)=G^{-1}(v(x))=G^{-1}(a\cdot x+b)$. Then $u$ is parallel. A priori, $\Omega$ could be a slab, but this is impossible since $\nabla u\neq 0$ by Step 1. This concludes the proof.
\end{proof}

\begin{remark}\label{re2.4}
During the proof Step 1, to obtain \cref{varphi} for $p>2$, we need to assume \cref{e1.1}. For $p\leq 2$ (including the Laplacian studied in \cite[Problem (1.1)]{Ruiz_Sicbaldi_Wu24}), this assumption is not needed since \cref{e1.1} can be derived directly from the fact that $F(u_0)=0$ and $u_0$ is the maximum point of $F$.
\end{remark}

Now, we can give the

\medskip

\noindent\textbf{Proof of \Cref{Th.P1}.} By \Cref{pr2.1}, we have \eqref{Modest}. In addition, from \Cref{esta}, we obtain the rest conclusion of \Cref{Th.P1}. ~\qed~\\

Next, we prove \Cref{Th.P2} and \Cref{coro}.

\medskip

\noindent\textbf{Proof of \Cref{Th.P2}.}
By the assumption \eqref{cond}, $Q$ attains its maximum at $\partial\Omega$. Then
    \[Q_{\nu}\geq0\quad \mbox{on  $\partial\Omega.$}\]
Let $H(q)$ be the mean curvature of $\partial\Omega$ at $q\in \partial \Omega$. By the boundary conditions $u=0$ and $u_{\nu}=-\kappa$ on $\partial \Omega$, we have
\begin{equation*}
u_i=-\kappa\nu_i~~\mbox{ on }~~\partial \Omega,~\forall ~1\leq i\leq n.
\end{equation*}
From the definition of $Q$ (see \cref{Pfunction}), we have
\begin{equation}\label{e.P_nu}
    \begin{aligned}
      Q_{\nu} &=(p-1)|\nabla u|^{p-2}u_ju_{ji} \nu_i+fu_i \nu_i\\
      &=-(p-1)\kappa^{p-1}u_{ij} \nu_i \nu_j-f\kappa \geq 0~~\mbox{ on}~~\partial \Omega.
    \end{aligned}
\end{equation}
Moreover, the solution $u$ satisfies
\begin{equation}\label{e.equation.boundary}
  \kappa^{p-2}\Delta u+(p-2) \kappa^{p-2} u_{ij} \nu_i \nu_j+f=0~~\mbox{ on}~~\partial \Omega.
\end{equation}
By using the fact (see \cite[Section 5.4]{MR615561} or \cite[Page 180]{MR4041100})
\begin{equation*}
  \frac{\partial^{2} u}{\partial\nu^{2}}=\Delta u-(n-1)H\frac{\partial u}{\partial\nu}
  ~~\mbox{ on }~~\partial\Omega,
\end{equation*}
we have
\begin{equation}\label{e.H2}
  \Delta u- u_{ij}\nu_i \nu_j +\kappa (n-1)H=0~~\mbox{ on}~~\partial \Omega.
\end{equation}
By considering \cref{e.equation.boundary}$-\kappa^{p-2}\times$\cref{e.H2}, we obtain
\begin{equation}\label{e.H3}
(p-1)\kappa^{p-2}u_{ij}\nu_i \nu_j- \kappa^{p-1}(n-1)H+f=0~~\mbox{ on}~~\partial \Omega.
\end{equation}
Finally, let \cref{e.P_nu} $+\kappa \times$ \cref{e.H3} and we have
\begin{equation*}
H\leq0~~\mbox{ on }~~\partial \Omega.
\end{equation*}

Since $\kappa\neq 0$, $u$ has no critical points close to $\partial \Omega$. Hence, $Q$ is a subsolution in a neighborhood of $\partial \Omega$. If $H(q)=0$ for some $q \in \partial \Omega$, then $Q_{\nu}(q)=0$. By the Hopf's lemma, $Q$ is a constant, at least in a neighborhood of $\partial \Omega$. We can now argue as in the proof of \Cref{esta} (see \cref{otra} and \cref{eq42}) to conclude that $u$ is parallel in such neighborhood. By unique continuation, $u$ is parallel and $\Omega$ is either a half-space or a slab. ~\qed~\\

\noindent\textbf{Proof of \Cref{coro}.} By \Cref{Th.P1}, we have
\begin{equation*}
  Q(x) \leq F(0)+\frac{p-1}{p}\kappa^p
\end{equation*}
for all $x \in \Omega$. Since $\kappa \neq 0$, according to \Cref{Th.P2}, $H \leq 0$. Moreover, if $H(q)=0$ at some point of $\partial \Omega$, then $Q$ is constant and $u$ is parallel. Additionally, by \Cref{esta}, we have $u$ has no critical points, and then $\Omega$ is a half-space. ~\qed~\\

\medskip

\noindent \textbf{Acknowledgements.}\\
Y. L. has been supported by the Grants PID2020-117868GB-I00 and PID2023-150727NB-I00 of the MICIN/AEI. J.W. has been supported by Proyecto de Consolidación Investigadora 2022, CNS2022-135640, MICINN.

\printbibliography

\end{document}